\documentclass[11pt,reqno]{amsart}
\usepackage{amssymb,amscd,amsbsy,mathrsfs}
\usepackage{amsthm}

\usepackage[russian]{babel}
\Russian
\usepackage[cp1251]{inputenc}
\usepackage[OT2,T1]{fontenc}

\newcommand{\lb}{\linebreak}

\newcommand{\e}{\varepsilon}

\newcommand{\z}{\zeta}

\newcommand{\s}{\sigma}
\renewcommand{\t}{\tau}

\renewcommand{\f}{\varphi}

\newcommand{\D}{\Delta}

\newcommand{\h}{{\mathscr H}}

\newcommand{\X}{{\mathscr X}}
\newcommand{\Y}{{\mathscr Y}}

\newcommand{\T}{{\Bbb T}}

\newcommand{\R}{{\Bbb R}}

\newcommand{\0}{{\boldsymbol{0}}}

\newcommand{\bs}{\boldsymbol}

\newcommand{\bS}{{\boldsymbol S}}

\newcommand{\rf}[1]{(\ref{#1})}

\newcommand{\df}{\stackrel{\mathrm{def}}{=}}

\newcommand{\supp}{\operatorname{supp}}

\newcommand{\trace}{\operatorname{trace}}

\newcommand{\const}{\operatorname{const}}

\newcommand{\eeq}{\end{equation}}
\newcommand{\beq}{\begin{equation}}
\newcommand{\bay}{\begin{eqnarray}}
\newcommand{\ba}{\begin{align*}}
\newcommand{\ea}{\end{align*}}
\newcommand{\ey}{\end{eqnarray}}
\newcommand{\bey}{\begin{eqnarray*}}
\newcommand{\eey}{\end{eqnarray*}}

\newcommand{\be}{\infty}

\newcommand{\bl}{\blacksquare}

\newcommand{\Pf}{{\bf Доказательство. }}

\newtheorem{thm}{\hspace{\parindent}Теорема}[section]

\newtheorem{lem}[thm]{\hspace{\parindent}Лемма}

\newcommand\OL{{\rm OL}}

\newcommand\fM{\frak M}

\newcommand\dg{\frak D}

\newcommand\mB{\mathcal{B}}

\begin{document}

\

\newcommand{\vse}{\vspace{.2in}}
\numberwithin{equation}{section}

\title{Формула следов Крейна для унитарных операторов и операторно липшицевы функции}
\author{А.Б. Александров и В.В.Пеллер}
\thanks{Исследования первого автора частично поддержаны грантом РФФИ 14-01-00198;
исследования второго автора частично поддержаны грантом NSF DMS 130092}

\begin{abstract}
Основной результат работы состоит в описании класса функций на единичной окружности, для которых справедлива формула следов Крейна для произвольных пар унитарных операторов с ядерной разностью. Этот класс состоит в точности из операторно липшицевых функция на окружности.
\end{abstract}

\maketitle

\tableofcontents

\

\setcounter{section}{0}
\section{\bf Введение}
\setcounter{equation}{0}
\label{Vved}

\

Функция спектрального сдвига для пары самосопряжённых операторов в гильбертовом пространстве была введена в работе И.М. Лифшица \cite{Li}. Там же была установлена формула следов для разности функций от возмущённого оператора и невозмущённого оператора. Идеи Лифшица были развиты в работе М.Г. Крейна \cite{Kr1}, в которой функция спектрального сдвига $\bs{\xi}$ класса $L^1(\R)$ определяется для произвольной пары самосопряжённых операторов $A$ и $B$ с ядерной разностью $A-B$ и доказывается формула следов 
\bay
\label{fsLK}
\trace\big(f(A)-f(B)\big)=\int_\R f'(t)\bs{\xi}(t)\,dt,
\ey
открытая Лифшицем, в значительно более общей ситуации, когда производная функции 
$f$ является преобразованием Фурье комплексной борелевской меры на $\R$.

Позже в работах \cite{Pe2} и \cite{Pe4} формула следов \rf{fsLK} была распространена на произвольные функции $f$ класса Бесова $B_{\be,1}^1(\R)$ (см. \cite{Pee} по поводу определения классов Бесова).

С другой стороны, очевидно, что правая часть формулы \rf{fsLK} имеет смысл для произвольных липшицевых функции $f$. М.Г. Крейн задал в работе \cite{Kr1} вопрос, можно ли обобщить формулу \rf{fsLK} на случай произвольных липшицевых функций $f$. Оказалось, что ответ на этот вопрос отрицателен: в работе \cite{F2} Ю.Б. Фарфоровская построила пример липшицевой функции $f$ и самосопряжённых операторов $A$ и $B$ таких, что $A-B$ входит в класс ядерных операторов $\bS_1$, но
$f(A)-f(B)\notin\bS_1$.

Таким образом, вопрос о применимости формулы следов \rf{fsLK} 
фактически распадается на два самостоятельных вопроса:

\medskip

(а) {\it Для каких функций $f$ на $\R$ выполняется условие
$$
A-B\in\bS_1\quad\Longrightarrow\quad f(A)-f(B)\in\bS_1
$$
для не обязательно ограниченных самосопряжённых операторов $A$ и $B$?}

(б) {\it Если $f$ -- функция, удовлетворяющая условию} (а){\it, то должны ли совпадать левая и правая части равенства {\em\rf{fsLK}}?}

Хорошо известно (см., например, недавний обзор \cite{AP}, теор. 3.6.5), что функция $f$ на $\R$ удовлетворяет условию (а) в том и только в том случае, когда она {\it операторно липшицева}, т.е.
имеет место неравенство
$$
\|f(A)-f(B)\|\le\const\|A-B\|
$$
для произвольных (ограниченных или неограниченных) самосопряжённых операторов $A$ и $B$.

То, что не всякая липшицева функция является операторно липшицевой, было обнаружено в работе Ю.Л. Фарфоровской \cite{F1}. Далее, в работе \cite{JW} было показано, что операторно липшицевы функции дифференцируемы в каждой точке. Это сразу же влечёт полученные ранее результаты работ \cite{Mc} и \cite{Ka}: функция $x\mapsto|x|$ не является операторной липшицевой. Отметим также, что операторно липшицевы функции не обязательно непрерывно дифференцируемы, что было показано в \cite{KS}.
В работах \cite{Pe2} и \cite{Pe4} были получены необходимые условия для операторной липшицевости. Эти необходимые условия основаны на описании ядерных операторов Ганкеля \cite{Pe1} (см. также \cite{Pe3}). 

Мы отсылаем читателя к недавнему обзору \cite{AP}, где проводится подробный анализ необходимых условий и достаточных условий для операторной липшицевости.

Положительный ответ на вопрос (б) был получен в недавней работе \cite{Pe5}: формула  \rf{fsLK} справедлива для произвольных операторно липшицевых функций $f$. Таким образом, класс функций, для которых справедлива формула следов \rf{fsLK} для любых самосопряжённых операторов $A$ и $B$ с ядерной разностью, совпадает с классом операторно липшицевых функций.

В этой работе мы предлагаем решение аналогичной задачи для функций от унитарных операторов.

Функция спектрального сдвига для пары унитарных операторов с ядерной разностью была определена в работе М.Г. Крейна \cite{Kr2} (см. также работу \cite{Kr3}, где приводится подробное изложение результатов). Пусть $U$ и $V$ -- унитарные операторы с ядерной разностью $U-V$. Тогда существует суммируемая функция $\bs{\xi}$ на единичной окружности $\T$ (называемая {\it функцией спектрального сдвига} для пары $(U,V)\:$) такая, что имеет место формула следов
\bay
\label{LKun}
\trace\big(f(U)-f(V)\big)=\int_\T f'(\z)\bs{\xi}(\z)\,d\z
\ey
для достаточно хороших функций $f$. В отличие от случая самосопряжённых операторов,
функция $\bs{\xi}$ не определяется однозначно парой $(U,V)$; она определяется с точностью до постоянной функции. Поэтому разумно потребовать, чтобы среднее значение функци $\bs{\xi}$ на окружность $\T$ было равно $0$.

В работе М.Г. Крейна \cite{Kr2} показано, что формула следов \rf{LKun} справедлива в случае, если производная $f'$ имеет абсолютно сходящийся ряд Фурье. В работе \cite{Pe2} формулу следов \rf{LKun} удалось распространить на функции $f$ класса Бесова $B_{\be,1}^1(\T)$.

Отметим, что так же, как и в случае функций от самосопряжённых операторов, функция $f$ обеспечивает ядерные приращения при ядерных возмущениях, т.е.
$$
U-V\in\bS_1\quad\Longrightarrow\quad f(U)-f(V)\in\bS_1
$$
в том и только в том случае, если $f$ -- {\it операторно липшицева функция}, т.е.
$$
\|f(U)-f(V)\|\le\const\|U-V\|
$$
для любых унитарных операторов $U$ и $V$. В классе $\OL(T)$ операторно липшицевых функций вводится естественная полунорма
$$
\|f\|_{\OL}\df\sup\frac{\|f(U)-f(V)\|}{\|U-V\|},
$$
где супремум берётся по всем унитарным операторам $U$ и $V$ таким, что $U\ne V$.

Основной результат этой работы будет получен в \S\:\ref{Unitar} и состоит в том, что формула следов \rf{LKun} справедлива для произвольных операторно липшицевых функций $f$. Ясно, что это -- максимальный класс функций с таким свойством. Из этого результата мы выведем следующий любопытный факт: функция
$$
\z\mapsto\trace\big(f(\z U)-f(\z V)\big)
$$
непрерывна на $\T$ для любой операторно липшицевой функции $f$ и для произвольной пары $(U,V)$ унитарных операторов с ядерной разностью.

Отметим, что доказательство, полученное в работе \cite{Pe5} для функций от самосопряжённых операторов, не распространяется на случай функций от унитарных операторов, ибо оно использует результат работы \cite{KPSS} о дифференцируемости операторных функций в норме Гильберта--Шмидта. Нам неизвестно, справедлив ли аналог этого утверждения в случае функций от унитарных операторов.

Вместо этого мы будем использовать в этой работе дифференцируемость соответствующих операторных функций в сильной операторной топологии, которая будет установлена в \S\:\ref{Difsil}.

В \S\:\ref{Dvoiopi} мы помещаем краткое введение двойных операторных интегралов. При этом мы приводим общую формулу следов, которая будет использоваться для доказательства основного результата в \S\:\ref{Unitar}.

Наконец, в \S\:\ref{Samosop} мы кратко обрисуем альтернативный подход в случае функций от самосопряжённых операторов, который в отличие от доказательства, полученного в работе \cite{Pe5}, вместо дифференцируемости в норме Гильберта--Шмидта использует дифференцируемость в сильной операторной топологии.

\

\section{\bf Двойные операторные интегралы и мультипликаторы Шура}
\setcounter{equation}{0}
\label{Dvoiopi}

\

Двойные операторные интегралы появились в работе Ю.Л. Далецкого и С.Г. Крейна \cite{DK}. Затем М.Ш. Бирман и М.З. Соломяк в работах \cite{BS1}, \cite{BS2} и \cite{BS4} создали красивую теорию двойных операторных интегралов.

Пусть $(\X,E_1)$ и $(\Y,E_2)$ -- пространства со спектральными мерами $E_1$ и $E_2$
в гильбертовом пространстве $\h$ и пусть $\Phi$ -- ограниченная измеримая функция на $\X\times\Y$. Двойные операторные интегралы -- это выражения вида
\bay
\label{doi}
\int\limits_\X\int\limits_\Y\Phi(x,y)\,d E_1(x)T\,dE_2(y).
\ey
Отправной точкой для определения двойных операторных интегралов в работах Бирмана и Соломяка являлся случай, когда  $T$ -- оператор Гильберта--Шмидта, и в этом случае двойные операторные интегралы определяются для произвольной измеримой ограниченной функции $\Phi$. 

Мы не будем здесь рассматривать случай операторов Гильберта--Шмидта, а отошлём читателя к обзору \cite{AP}, глава II,
в котором двойные операторные интегралы подробно обсуждены. 

Для того, чтобы определить двойные операторные интегралы вида \rf{doi} для произвольных ограниченных операторов $T$, нужно наложить на функцию $\Phi$ дополнительные ограничения. Именно, двойные операторные интегралы вида \rf{doi}
могут быть определены для любых ограниченных операторов $T$ при условии, что $\Phi$ входит в класс {\it мультипликаторов Шура} $\fM(E_1,E_2)$ по  отношению к спектральным мерам $E_1$ и $E_2$. Класс $\fM(E_1,E_2)$ допускает различные описания, см. \cite{Pe2}, \cite{Pi} и \cite{AP}. 

Здесь мы приведём одно такое описание: {\it $\Phi\in\fM(E_1,E_2)$ в том и только в том случае, когда $\Phi$ принадлежит тензорному произведению Хогерупа 
$L^\be(E_1)\otimes_{\rm h}L^\be(E_2)$, т.е. $\Phi$ допускает представление
\bay
\label{predstava}
\Phi(x,y)=\sum_n\f_n(x)\psi_n(y),
\ey
где $\f_n$ и $\psi_n$ -- измеримые функции, удовлетворяющие условию}
$$
\Big\|\sum_n|\f_n|^2\Big\|_{L^\be(E_1)}\le\|\Phi\|_{\fM(E_1,E_2)}
\quad\mbox{и}\quad
\Big\|\sum_n|\psi_n|^2\Big\|_{L^\be(E_2)}\le\|\Phi\|_{\fM(E_1,E_2)},
$$
где $\|\Phi\|_{\fM(E_1,E_2)}$ -- норма трансформатора
$$
T\mapsto\iint\Phi\,dE_1T\,dE_2
$$
в пространстве операторов в гильбертовом пространстве.
При этом имеет место равенство
\bay
\label{doitproh}
\int\limits_\X\int\limits_\Y\Phi(x,y)\,d E_1(x)T\,dE_2(y)=
\sum_n\Big(\int\f_n\,dE_1\Big)T\Big(\int\psi_n\,dE_2\Big),
\ey
причём ряд в правой части равенства сходится в слабой операторной топологии, и его сумма не зависит от выбора представления \rf{predstava}.

Если $\Phi\in\fM(E_1,E_2)$, а $T$ -- ядерный оператор, то двойной операторный
интеграл \rf{doi} также должен быть ядерным и при этом справедливо неравенство
\bay
\label{yadotse}
\left\|\,\int\limits_\X\int\limits_\Y\Phi(x,y)\,d E_1(x)T\,dE_2(y)\right\|_{\bS_1}
\le\|\Phi\|_{\fM(E_1,E_2)}\|T\|_{\bS_1}.
\ey


Предположим, что $f$ -- операторно липшицева функция на единичной окружности $\T$.
Рассмотрим её разделённую разность на $\T\times\T$:
$$
(\dg f)(\z,\t)\df\left\{\begin{array}{ll}\frac{f(\z)-f(\t)}{\z-\t},&\z\ne\t,\\[.2cm]
f'(\z),&\z=\t
\end{array}\right.
$$
(в силу результатов работы \cite{JW} операторно липшицевы функции на окружности дифференцируемы в каждой точке). Хорошо известно, что в этом случае разделённая разность $\dg f$ является мультипликатором Шура для любых борелевских спектральных мер $E_1$ и $E_2$, причём верно и обратное утверждение: если функция $f$ на $\T$ всюду дифференцируема и $\dg f$ -- мультипликатор Шура для любых борелевских спектральных мер, то функция $f$ операторно липшицева (см., например, обзор
\cite{AP}, теор. 3.3.6). Более того, имеет место равенство
$$
\|f\|_{\OL}=\sup\|\dg f\|_{\fM(E_1,E_2)},
$$
где супремум берётся по всем борелевским спектральным мерам $E_1$ и $E_2$ на $\T$.

Также хорошо известно (cм. \cite{BS4} и обзор \cite{AP}), что при этих условиях имеет место формула
\bay
\label{DKBS}
f(U)-f(V)=\iint\limits_{\T\times\T}(\dg f)(\z,\t)\,dE_U(\z)(U-V)\,dE_V(\T),
\ey
где $E_U$ и $E_V$ -- спектральные меры операторов $U$ и $V$.

Предположим теперь, что $E_1$ и $E_2$ -- борелевские спектральные меры на локально компактных топологических пространствах $\X$ и $\Y$, по крайней мере одно из которых сепарабельно и пусть $\supp E_1=\X$, а $\supp E_2=\Y$. Тогда, если принять во внимание теорему 2.1 работы \cite{AP2}, из теоремы 2.2.4 работы \cite{AP} получаем следующее утверждение:

{\it Пусть $\Phi$ -- функция на $\X\times\Y$, непрерывная по каждой переменной.
Тогда \lb$\Phi\in\fM(E_1,E_2)$ в том и только в том случае, когда она входит в тензорное произведение Хогерупа 
$C_{\rm b}(\X)\!\otimes_{\rm h}\!C_{\rm b}(\Y)$ пространств 
$C_{\rm b}(\X)$ и $C_{\rm b}(\Y)$
ограниченных непрерывных функций на $\X$ и $\Y$, т.е. $\Phi$ допускает представление вида 
{\em\rf{predstava}}, в котором
$\f_n\in C_{\rm b}(\X)$, $\psi_n\in C_{\rm b}(\Y)$ и имеют место неравенства}
$$
\sum_n|\f_n|^2\le\|\Phi\|_{\fM(E_1,E_2)}
\quad\mbox{и}\quad\sum_n|\psi_n|^2\le\|\Phi\|_{\fM(E_1,E_2)}.
$$

Перейдём теперь к общей формуле следов для двойных операторных интегралов.

Пусть $T$ -- ядерный оператор в гильбертовом пространстве, а
$E$ -- спектральная мера на $\s$-алгебре подмножеств множества $\X$, а
$\Phi\in\fM(E,E)$.
Вычислим след двойного операторного интеграла
$$
\iint\Phi(x,y)\,dE(x)T\,dE(y).
$$
В работе \cite{BS4} найдена следующая формула:
\bay
\label{BStf}
\trace\left(\iint\Phi(x,y)\,dE(x)T\,dE(y)\right)=\int\Phi(x,x)\,d\mu(x),
\ey
где $\mu$ -- комплексная мера на этой же $\s$-алгебре, определённая равенством
$$
\mu(\D)=\trace\big(TE(\D)\big).
$$

Для обоснования правой части формулы \rf{BStf} нужно понять, как можно интерпретировать
значения функции $\Phi$ на диагонали $\{(x,x):~x\in\X\}$.
В работе
\cite{Pe+} была дана следующая интерпретация формулы \rf{BStf}. Мы можем определить след ${\mathscr T}\Phi$ функции $\Phi$ класса $\fM(E,E)$  на диагонали равенством
$$
({\mathscr T}\Phi)(x)\df\sum_n\f_n(x)\psi_n(x),
$$
где $\f_n$ и $\psi_n$ функции из представления \rf{predstava} функции $\Phi$ в виде 
тензорного произведения Хогерупа 
$L^\be(E)\otimes_{\rm h}L^\be(E)$.
Тогда след ${\mathscr T}\Phi$ функции $\Phi$ класса $\fM(E,E)$ на диагонали принадлежит пространству $L^\be(E)$ и не зависит от выбора представления 
\rf{predstava}.
В формуле \rf{BStf} следует понимать $\Phi(x,x)$, как $({\mathscr T}\Phi)(x)$,
см. \cite{Pe5}, \S\:1.1.

Наконец, предположим, что $E$ -- борелевская спектральная мера на локально компактном топологическом пространстве $\X$, а  
$\Phi$ -- функция на $\X\times\X$, непрерывная по каждой переменной.  Тогда справедливо следующее утверждение (см. \cite{Pe5}):

\begin{thm}
\label{sleddvoi}
Пусть  $E$ -- борелевская спектральная мера на локально компактном пространстве $\X$, а 
$\Phi$ -- функция класса $\fM(E,E)$. Если функция $\Phi$ непрерывна по каждой переменной, то для любого ядерного оператора $T$ справедлива
формула {\em\rf{BStf}}.
\end{thm}

Действительно, достаточно рассмотреть случай, когда $\supp E=\X$ и рассмотреть
представление \rf{predstava} функции $\Phi$ в виде тензорного произведения Хогерупа 
\lb$C_{\rm b}(\X)\!\otimes_{\rm h}\!C_{\rm b}(\X)$. Легко видеть, что в этом случае 
$({\mathscr T}\Phi)(x)=\Phi(x,x)$, $x\in\X$.

\

\section{\bf Операторная дифференцируемость в сильной операторной топологии}
\setcounter{equation}{0}
\label{Difsil}

\

В этом параграфе для операторно липшицевой функций $f$ на $\T$, унитарного оператора $U$ и ограниченного самосопряжённого оператора $A$ мы рассмотрим
задачу дифференцируемости операторной функции
$$
t\mapsto f\big(e^{{\rm i}tA}U\big)
$$
в сильной операторной топологии. Отметим здесь, что аналог следующей теоремы для функций от самосопряжённых операторов был установлен в \cite{AP}, теор. 3.5.5; см. также \S\:\ref{Samosop} этой статьи.

\begin{thm} 
\label{sildiftorn}
Пусть $f$ -- операторно липшицева функция на $\T$,  $U$ -- унитарный оператор, а
$A$  --  ограниченный самосопряжённый оператор.
Тогда
\bay
\label{vsilopto}
\lim_{t\to0}\frac1t\Big(f\big(e^{{\rm i}tA}U\big)-f(U)\Big)={\rm i}\int_\T\int_\T\t(\dg f)(\z,\t)\,dE_U(\z)A\,dE_U(\t),
\ey
где предел берётся в сильной операторной топологии.
\end{thm}

Заметим, что в работе \cite{Pe2} формула \ref{vsilopto} получена для функций $f$ класса Бесова $B_{\be,1}^1(\T)$, при этом предел в левой части равенства \ref{vsilopto} существует по операторной норме.

Нам понадобится следующее вспомогательное утверждение.

\begin{lem} 
\label{Xnun-n}
Пусть $\{X_n\}_{n\ge0}$ -- последовательность в пространстве $\mB(\h)$ линейных ограниченных операторов в гильбертовом пространстве $\h$,
а $\{u_n\}_{n\ge0}$ -- последовательность в $\h$. Предположим, что 
$$
\sum_{n\ge0}X_nX_n^*\le a^2I\quad\mbox{и}\quad
 \sum_{n\ge0}\|u_n\|^2\le b^2
$$
для неотрицательных чисел $a$ и $b$. Тогда ряд $\sum_{n\ge0}X_nu_n$
слабо сходится и 
$$
\Big\|\sum_{n\ge0}X_nu_n\Big\|\le ab.
$$
\end{lem}

\Pf Пусть $v\in\h$ и $\|v\|=1$. Тогда
$$
\sum_{n\ge0}|(X_nu_n,v)|=\sum_{n\ge0}|(u_n,X_n^*v)|\le
\Big(\sum_{n\ge0}\|u_n\|^2\Big)^{1/2}\Big(\sum_{n\ge0}\|X_n^*v\|^2\Big)^{1/2}\le ab,
$$
откуда следует доказываемое утверждение. $\bl$

\medskip

{\bf Доказательство теоремы \ref{sildiftorn}.} Как  мы упоминали в 
\S\:\ref{Dvoiopi}, функция $f$ дифференцируема всюду на $\T$, разделённая разность 
$\dg f$ является мультипликатором Шура для любых борелевских спектральных мер $E_1$ и $E_2$, причём
$$
\|\dg f\|_{\fM(E_1,E_2)}\le\|f\|_{\OL}.
$$
Мы также отмечали в \S\:\ref{Dvoiopi}, что 
существуют последовательности непрерывных на 
$\T$
функций $\{\f_n\}_{n\ge0}$ и $\{\psi_n\}_{n\ge0}$ таких, что 

{\rm a)} $\sum\limits_{n\ge0}|\f_n|^2\le\|f\|_{\OL(\T)}$ всюду на $\T$,

{\rm б)} $\sum\limits_{n\ge0}|\psi_n|^2\le\|f\|_{\OL(\T)}$ всюду на $\T$,

{\rm в)} $(\dg f)(\z,\t)=\sum\limits_{n\ge0}\f_n(\z)\psi_n(\t)$ при всех $\z$ и $\t$ из $\T$.

Ввиду тождеств \rf{doitproh} и \rf{DKBS} мы должны показать, что
$$
\lim_{t\to0}\frac1t\sum_{n\ge0}\f_n(e^{{\rm i}tA}U)(e^{{\rm i}tA}-I)U\psi_n(U)={\rm i}\sum_{n\ge0}\f_n(U)AU\psi_n(U).
$$
Причём ряды суммируются в слабой операторной топологии, а предел берётся
в сильной операторной топологии. Заметим, что $\lim\limits_{t\to0}t^{-1}(e^{{\rm i}tA}-I)={\rm i}A$ 
по операторной норме. Таким образом, достаточно доказать, что
$$
\lim_{t\to0}\sum_{n\ge0}\f_n(e^{{\rm i}tA}U)AU\psi_n(U)=\sum_{n\ge0}\f_n(U)AU\psi_n(U)
$$
в сильной операторной топологии.
Иными словами, нам нужно доказать, что для любого вектора $u\in\h$ мы имеем
$$
\lim_{t\to0}\sum_{n\ge0}(\f_n(e^{{\rm i}tA}U)-\f_n(U))AU\psi_n(U)u=\0,
$$
где ряд суммируется в слабой топологии пространства $\h$, а предел берётся в пространстве $\h$ по норме.
Будем считать, что $\|u\|=1$ и $\|f\|_{\OL(\T)}=1$.  Тогда $\sum_{n\ge0}|\f_n|^2\le1$ и
$\sum_{n\ge0}|\psi_n|^2\le1$ всюду на $\T$.

Положим $u_n\df AU\psi_n(U)u$. Имеем:
$$
\sum_{n\ge0}\|u_n\|^2\le\|A\|^2\sum_{n\ge0}\|\psi_n(U)u\|^2=
\|A\|^2\sum_{n\ge0}(|\psi_n|^2(U)u,u)\le\|A\|^2<+\be.
$$
Пусть $\e>0$. Выберем  натуральное число $N$ так, чтобы 
$\sum_{n>N}\|u_n\|^2<\e^2$. Тогда из леммы  \ref{Xnun-n}
следует, что
$$
\Big\|\sum_{n>N}(\f_n(e^{{\rm i}tA}U)-\f_n(U))u_n\Big\|\le2\e
$$
при всех $t\in\R$. 

Легко видеть и хорошо известно, что если $h$ -- непрерывная функция на $\T$, то отображение
$$
U\mapsto h(U)
$$
непрерывно на множестве унитарных операторов в операторной норме (достаточно аппроксимировать функцию $h$ тригонометрическими полиномами).

Тогда
$$
\left\|\sum_{n=0}^N(\f_n(e^{{\rm i}tA}U)-\f_n(U))u_n\right\|
\le\|A\|\sum_{n=0}^N\Big\|\f_n\big(e^{{\rm i}tA}U\big)-\f_n(U)\Big\|<\e
$$
при всех достаточно близких к нулю $t$. Таким образом, 
$$
\Big\|\sum_{n\ge0}\big(\f_n(e^{{\rm i}tA}U)-\f_n(U)\big)u_n\Big\|<3\e
$$
при всех достаточно близких к нулю $t$. $\bl$

\

\section{\bf Формула следов и операторная липшицевость}
\setcounter{equation}{0}
\label{Unitar}

\

В этом параграфе мы установим основной результат работы, которые состоит в следующем:

\begin{thm}
\label{osnrezu}
Формула следов {\em\rf{LKun}} справедлива для любой операторно липшицевой функции $f$ на $\T$ и для любой пары $(U,V)$ унитарных операторов с ядерной разностью
$U-V$.
\end{thm}

В доказательстве будет использоваться идея Бирмана и Соломяка из работы \cite{BS3}, используемая ими для их подхода к построению функции спектрального сдвига. При этом в их работе накладываются более обременительные ограничения на функцию $f$.

\medskip

{\bf Доказательство теоремы \ref{osnrezu}.} Прежде всего, легко видеть, что при условии $U-V\in\bS_1$ существует ядерный самосопряжённый оператор $A$ такой, что
$V=e^{{\rm i}A}U$. В силу теоремы \ref{sildiftorn} функция 
$t\mapsto f\big(e^{{\rm i}tA}U\big)$ дифференцируема в сильной операторной топологии и
\bay
\label{formdlyapro}
Q_s\df\frac{d}{dt}f\big(e^{{\rm i}tA}U\big)\Big|_{t=s}=
{\rm i}\int_\T\int_\T\t(\dg f)(\z,\t)\,dE_s(\z)A\,dE_s(\t),
\ey
где $E_s$ -- спектральная мера унитарного оператора $V_s\df e^{{\rm i}sA}U$.

Как мы отмечали в \S\:\ref{Dvoiopi}, разделённая разность $\dg f$ является мультипликатором Шура, а посему, ввиду \rf{yadotse},
$$
Q_s\in\bS_1\quad\mbox{и}\quad\sup_{s\in[0,1]}\|Q_s\|_{\bS_1}<\be.
$$

Из определения функции $s\mapsto Q_s$ следует измеримость функции $s\mapsto Q_su$ для любого
вектора $u$ из гильбертова пространства $\h$. Следовательно, скалярная функция $s\mapsto (Q_su,v)$
измерима для любых $u$ и $v$ из $\h$. Отсюда легко вытекает измеримость функции $s\mapsto \trace(Q_sT)$
для любого оператора $T\in\mB(\h)$. Таким образом, $\bS_1$-значная функция $s\mapsto Q_s$
слабо измерима, а тогда она измерима и в сильном смысле, поскольку пространство $\bS_1$  сепарабельно, см., например, \cite{I}, глава V, \S\:4.

Теперь из равенства \rf{formdlyapro} вытекает, что
$$
f(V)-f(U)=\int_0^1Q_s\,ds,
$$
где интеграл понимается, как интеграл Бохнера в пространстве $\bS_1$.

Тогда
$$
\trace\big(f(V)-f(U)\big)=\int_0^1\trace Q_s\,ds.
$$ 

По теореме \ref{sleddvoi}
$$
\trace Q_s=\int_\T \z f'(\z)\,d\nu_s(\z),
$$
где $\nu_s$ -- комплексная борелевская мера на $\T$, определённая равенством
$$
\nu_s(\D)\df\trace(E_s(\D)A)
$$
для борелевского подмножества $\D$ окружности $\T$.

Мы можем отождествить пространство ${\mathcal M}(\T)$ комплексных борелевских мер на $\T$ с двойственным пространством к пространству $C(\T)$ непрерывных функций на $\T$. Покажем, что функция $s\mapsto\nu_s$ непрерывна в слабой топологии 
$\s\big({\mathcal M}(\T),C(\T)\big)$. Действительно, если $h\in C(\T)$, то
$$
\int_\T h\,d\nu_s=\trace(h(V_s)A).
$$
Как мы уже отмечали в доказательстве теоремы \ref{sildiftorn}, функция $s\mapsto h(V_s)$ непрерывна в операторной норме. Отсюда и вытекает слабая непрерывность
непрерывность функции $s\mapsto\nu_s$.

Определим теперь комплексную меру $\nu$ равенством
$$
\nu=-\int_0^1\nu_s\,ds.
$$
Здесь интеграл понимается, как интеграл непрерывной функции в топологии 
\lb$\s\big({\mathcal M}(\T),C(\T)\big)$.

Тогда
$$
\trace\big(f(U)-f(V)\big)=\int_\T\z f'(\z)\,d\nu(\z).
$$
С другой стороны, для тригонометрических полиномов $f$ имеет место равенство
$$
\trace\big(f(U)-f(V)\big)=\int_\T f'(\z)\bs{\xi(\z)}\,d\z.
$$
Отсюда следует, что существует константа $c$ такая, что
$$
\z\,d\nu(\z)=\bs{\xi}(\z)\,d\z+c\,\z^{-1}\,d\z,
$$
что завершает доказательство теоремы. $\bl$

\medskip

Теорема \ref{osnrezu} позволяет нам получить следующий занятный факт:

\begin{thm}
Пусть $f$ -- операторно липшицева функция на $\T$, а $U$ и $V$ -- унитарные операторы такие, что $U-V\in\bS_1$. Тогда функция
$$
\z\mapsto\trace\big(f(\z U)-f(\z V)\big),\quad\z\in\T,
$$
является непрерывной на $\T$.
\end{thm}

\Pf Пусть $f\in\OL(\T)$. Тогда для $\z\in\T$ положим $f_\z(\t)=f(\z\t)$, $\t\in\T$.
Имеем
$$
\trace\big(f(\z U)-f(\z V)\big)=\trace\big(f_\z(U)-f_\z(V)\big)
=\int_\T f_\z'(\t)\bs{\xi}(\t)\,d\t,
$$
где $\bs{\xi}$ -- функция спектрального сдвига для пары $(U,V)$. Осталось заметить, что функция
$$
\z\mapsto\int_\T f_\z'(\t)\bs{\xi}(\t)\,d\t,\quad\z\in\T,
$$
непрерывна на $\T$, ибо функция $\bs{\xi}$ суммируема, а функция $f'$ входит в $L^\be$. $\bl$

\

\section{\bf Альтернативный подход к случаю самосопряжённых операторов}
\setcounter{equation}{0}
\label{Samosop}

\

В работе \cite{Pe5} для того, чтобы доказать формулу следов Лифшица--Крейна для операторно липшицевых функций от самосопряжённых операторов, использовался следующий результат работы \cite{KPSS}: пусть $A$ -- самосопряжённый оператор, а $K$ -- самосопряжённый оператор класса $\bS_2$. Тогда, если $f$ -- всюду дифференцируемая функция на $\R$ с ограниченной производной, то функция
$t\mapsto f(A+tK)-f(A)$ дифференцируема по норме $\bS_2$ и
\bay
\label{kommunisticheskaya}
\frac d{dt}\big(f(A+tK)-f(A)\big)\Big|_{t=0}=
\int_\R\int_\R(\dg f)(x,y)\,dE_A(x)K\,dE_A(y).
\ey

Аналогом этого утверждения для функций от унитарных операторов было бы следующее утверждение: если $f$ -- всюду дифференцируемая функция на $\T$ с ограниченной производной, $U$ -- унитарный оператор, а $A$ -- самосопряжённый оператор класса $\bS_2$, то функция $t\mapsto f\big(e^{{\rm i}tA}\big)U$ дифференцируема по норме $\bS_2$ и имеет место равенство \rf{vsilopto}. К сожалению, мы не знаем, верно ли это утверждение.

Вместо этого утверждение мы использовали в этой работе дифференцируемость этой функции в сильной операторной топологии в случае, когда $f$ -- операторно липшицева функция на $\T$, см. теорему \ref{sildiftorn}.

В работе \cite{AP} (см. теор. 3.5.6) был получен следующий аналог теоремы \ref{sildiftorn} для функций от самосопряжённых операторов:

{\it Пусть $f$ -- операторно липшицева функция на $\R$, а $A$ и $K$ -- самосопряжённые операторы, причём оператор $K$ ограничен. Тогда функция
$t\mapsto\big(f(A+tK)-f(A)\big)$ дифференцируема в сильной операторной топологии,
и имеет место формула} \rf{kommunisticheskaya}.

\medskip

Эта теорема позволяет получить новое доказательство формулы следов Лифшица--Крейна для операторно липшицевых функций от самосопряжённых операторов, которое не использует упомянутый выше результат работы \cite{KPSS} о дифференцируемости операторных функций в норме Гильберта--Шмидта.

\

\noindent
\begin{tabular}{p{9cm}p{15cm}}
А.Б. Александров & В.В. Пеллер \\
Санкт-Петербургское отделение & Department of Mathematics \\
Математический институт Стеклова РАН  & Michigan State University \\
Фонтанка 27, 191023 Санкт-Петербург & East Lansing, Michigan 48824\\
Россия&USA
\end{tabular}


\begin{thebibliography}{99}
\label{bibl}







\bibitem{AP}{\sc А.Б. Александров и В.В. Пеллер}, {\em Операторно липшицевы функции}, Успехи Матем. Наук.






\bibitem{BS1} {\sc М.Ш. Бирман} и {\sc М.З. Соломяк},
{\em Двойные операторные интегралы Стилтьеса}, Проблемы
мат. физики. {\bf1}. Спектральная теория и волновые процесcы. Издат. ЛГУ (1966), 33 -- 67.

\bibitem{BS2} {\sc М.Ш. Бирман} и {\sc М.З. Соломяк},
{\em Двойные операторные интегралы Стилтьеса. II}, Проблемы
мат. физики. {\bf2}. Спектральная теория, проблемы дифракции. Издат. ЛГУ (1967), 26 -- 60.

\bibitem{BS3} {\sc М.Ш. Бирман} и {\sc М.З. Соломяк}, {\em Замечания о функции спектрального сдвига}, Записки Научн. Семин. ЛОМИ {\bf27} (1972), 33--46.

\bibitem{BS4} {\sc М.Ш. Бирман} и {\sc М.З. Соломяк},
{\em Двойные операторные интегралы Стилтьеса. III. Предельный переход под знаком интеграла}, 
Проблемы мат. физики. {\bf6}. Теория функций. Спектральная теория. 
Распространение волн.  Издат. ЛГУ (1973),  27--53.

\bibitem{DK} {\sc Ю.Л, Далецкий} и {\sc С.Г. Крейн}, {\it Интегрирование и дифференцирование эрмитовых операторов и приложение к теории возмущений}, Труды семинара по функц. анализу, Воронеж, 1956, т. 1, с. 81--106.




\bibitem{I} {\sc К.~Иосида},
{\em Функциональный анализ}.
Мир, М., 1967.







\bibitem{Kr1} {\sc М.Г. Крейн}, {\em О формуле следов в теории возмущений},
Мат. Сборник {\bf33} (1953), 597--626.

\bibitem{Kr2} {\sc М.Г. Крейн}, {\em Об определителях возмущения и формуле следов для унитарных и самосопряжённых операторов}, Док. АН СССР {\bf144:2} (1962),  268--271.

\bibitem{Kr3} {\sc М.Г. Крейн}, {\em О некоторых новых исследованиях по теории возмущений самосопряжённых операторов}. В книге: Первая летняя математическая школа, Киев, 1964, 103--187.

\bibitem{Li}{\sc И.М. Лифшиц}, {\em Об одной задаче теории возмущений, связанной с квантовой статистикой}, УМН {\bf 7:1(47)} (1952), 171--180.


\bibitem{Pe1} {\sc В.В. Пеллер}, {\em Операторы Ганкеля класса $\mathfrak S_p$ и
их приложения (рациональная аппроксимация, гауссовские процессы, проблема мажорации
операторов)}, Матем. сб.,
{\bf 113(155):4(12)} (1980), 538-581.


\bibitem{Pe2} {\sc В.В. Пеллер}, {\em Операторы Ганкеля в теории возмущений
унитарных и самосопряженных операторов}, Функц. анал. и его прил. {\bf19:2}  (1985),
37 -- 51.

\bibitem{Pe3} {\sc В.В. Пеллер}, {\em Операторы Ганкеля и их приложения},
НИЦ ``Регулярная и хаотическая динамика'', Институт компьютерных исследований, Москва,
Ижевск, 2005.



%

%
%

\bibitem{F1} {\sc Ю.Б. Фарфоровская}, {\em О связи метрики Канторовича-Рубинштейна для спектральных разложений самосопряженных операторов с функциями от операторов},
Вестник ЛГУ, {\bf 19} (1968), 94--97.

\bibitem{F2} {\sc Ю.Б. Фарфоровская}, {\em  Пример липшицевой функции от самосопряженного оператора, дающей неядерные приращения при ядерном возмущении},  Зап. научн. сем. ЛОМИ, {\bf30}  (1972), 146--153.



%
%
%

\bibitem{AP2}  {\sc A.B. Aleksandrov} and {\sc V.V. Peller},
{\em Estimates of operator moduli of continuity}.  J. Funct. Anal. {\bf261:10} (2011), 2741 -- 2796.




\bibitem{JW} {\sc B.E. Johnson} and {\sc J.P. Williams}, {\em The range of a normal derivation}.
Pacific J. Math. {\bf58} (1975), 105 -- 122.


%


\bibitem{Ka}  {\sc T. Kato}, {\em Continuity of the map $S\mapsto |S| $ for linear operators,}  Proc. Japan Acad.  {\bf49}  (1973), 157 -- 160.

\bibitem{KPSS} {\sc E. Kissin, D. Potapov, V. S. Shulman} and {\sc F. Sukochev}, 
{\it Operator smoothness in Schatten norms for functions of several variables: 
Lipschitz conditions, differentiability and unbounded derivations},  
Proc. Lond. Math. Soc. (3) {\bf105} (2012), 661--702.


\bibitem{KS} {\sc E. Kissin} and {V.S. Shulman}, {\em On a problem of
J. P. Williams}. Proc. Amer. Math. Soc. {\bf130} (2002), 3605 -- 3608.


\bibitem{Mc}
{\sc A. McIntosh}, {\em Counterexample to a question on commutators},
Proc. Amer. Math. Soc. {\bf29} (1971), 337 -- 340.



\bibitem{Pee} {\sc J. Peetre},
{\em New thoughts on Besov spaces}, Duke Univ. Press., Durham, NC, 1976.


\bibitem{Pe4} {\sc V.V. Peller} {\em Hankel operators in the perturbation theory of unbounded self-adjoint operators}.  Analysis and partial differential equations,  529 -- 544,
Lecture Notes in Pure and Appl. Math., {\bf122}, Dekker, New York, 1990.

%



\bibitem{Pe+} {\sc V.V. Peller}, {\em Multiple operator integrals in perturbation theory}, Bull. Math. Sci. {\bf6} (2016), 15--88.

\bibitem{Pe5} {\sc V.V. Peller},
{\em The Lifshits--Krein trace formula and operator Lipschitz functions},
Proc. Amer. Math. Soc.

\bibitem{Pi} {\sc G. Pisier}, {\em Similarity problems and completely bounded maps},
Second, expanded edition. Includes the solution to ``The Halmos problem''. Lecture Notes in Mathematics,
1618. Springer-Verlag, Berlin, 2001.




\end{thebibliography}
\end{document}